\documentclass[12pt]{article}
\usepackage{amsfonts}
\usepackage{amsmath}
\usepackage{latexsym}
\usepackage{amscd}

\newtheorem{theorem}{Theorem}

\newcommand{\be}{\begin{equation}}
\newcommand{\ee}{\end{equation}}

\begin{document}
\bibliographystyle{plain}

\thispagestyle{empty}
\setcounter{page}{0}

\vspace {2cm}

{\Large 
G. Morvai and B. Weiss:  Limitations on intermittent forecasting.}

\vspace {2cm}

{\Large Appeared in :   Statist. Probab. Lett.  72  (2005),  no. 4, pp. 285--290. }

\vspace {2cm}

\begin{abstract}
Bailey showed that the general pointwise forecasting for stationary and ergodic 
time series has a 
negative solution.  However, it is  known that for Markov chains the problem can be solved. 
 Morvai showed  that there is a stopping time sequence $\{\lambda_n\}$ such that 
 $P(X_{\lambda_n+1}=1|X_0,\dots,X_{\lambda_n} ) $ can be estimated from samples $(X_0,\dots,X_{\lambda_n})$ such that 
 the difference between the conditional probability and the estimate  vanishes 
 along these stoppping times for all stationary and ergodic binary time series. 
We will show it is not possible to estimate 
the above conditional probability
along a stopping time sequence 
  for all stationary and ergodic binary time series in a pointwise sense such that 
 if the time series turns out to be a Markov chain, 
the predictor will predict eventually for all $n$. 
\end{abstract}



\pagebreak

\section{Introduction and Statement of  Results}

Cover \cite{Cover75} posed the following fundamental problem  concerning 
forecasting for stationary and ergodic binary  time series $\{ X_n\}_{n=-\infty}^{\infty}$.
(Note that a stationary time series $\{ X_n\}_{n=0}^{\infty}$ 
can be extended to be a two sided stationary time series 
$\{ X_n\}_{n=-\infty}^{\infty}$.)

\bigskip
\noindent
{\bf Problem 1} 

{\it Is there an estimation scheme $f_{n}$ for the value 
$P(X_{n+1}=1|X_0,X_1,\dots,X_n)$ such that 
$f_{n}$ depends solely on the data segment $(X_0,X_1,\dots, X_n)$ and 
$$
\lim_{n\to\infty} |f_{n}(X_0,X_1,\dots,X_n)-P(X_{n+1}=1|X_0,X_1,\dots,X_n)|=0
$$
almost surely for all stationary and ergodic binary    
time series $\{ X_n\}_{n=-\infty}^{\infty}$? }

\bigskip
\noindent
This problem was answered by Bailey \cite{Bailey76} in a negative way, that is, he showed that there is no such 
scheme. 
(Also see Ryabko \cite{Ryabko88}, Gy\"orfi,  Morvai, Yakowitz \cite{GYMY98}
 and  Weiss \cite{Weiss00}.)

\bigskip
\noindent
Morvai \cite{Mo00} considered the following modification of Problem ~1.

\bigskip 
\noindent
{\bf Problem 2} 

{\it 
 Are  there a strictly  increasing sequence of stopping times $\{\lambda_n\}$ and estimators  
$\{h_n(X_{0}, \dots, X_{\lambda_n})\}$ such that  
for all  stationary ergodic  binary time series $\{X_n\}$ 
 the estimator $h_n$   is consistent at stopping times $\lambda_n$, 
that is,   
$$
\lim_{n\to\infty} | h_n(X_0, \dots, X_{\lambda_n})- 
P(X_{\lambda_n+1}=1|X_0,\dots,X_{\lambda_n})|=0  
$$
almost surely ?  }

\bigskip
\noindent
Morvai \cite{Mo00} constructed a scheme that solves Problem~2. Unfortunatelly, 
his  stopping times  grow
extremly rapidly and so that scheme is not practical at all. 

\bigskip
\noindent
Let ${\cal X}^{*-}$ be the set of all one-sided  binary sequences, that is, 
$$
{\cal X}^{*-}=\{ (\dots,x_{-1},x_0): x_i\in \{0,1\} \ \ \mbox{for all $-\infty<i\le 0$}\}.
$$

\noindent
Define the distance $d^*(\cdot,\cdot)$ on ${\cal X}^{*-}$  as follows. Let 
$$
d^*((\dots,x_{-1},x_{0}),(\dots,y_{-1},y_{0}))=
\sum_{i=0}^{\infty} 2^{-i-1} |x_{-i}-y_{-i}|.
$$

\bigskip
\noindent
{\sc Definition}
The conditional probability $P(X_1=1|\dots,X_{-1},X_0 )$
is almost surely continuous if to some set $C\subseteq {\cal X}^{*-}$  
which has probability one   
the conditional probability $P(X_1=1|\dots,X_{-1},X_0 )$ restricted to this set $C$ 
is continuous with respect to metric $ d^*(\cdot,\cdot)$. 

\bigskip
The processes with almost surely continuous conditional probability 
generalizes the processes for which it is actually continuous, these are 
 essentially the Random Markov Processes of Kalikow~ 
\cite{Ka90}, or the continuous g-measures studied by Mike Keane  \cite{Ke72}.

\bigskip
A more moderate growth ( compared to Morvai \cite{Mo00} ) was achieved by Morvai and Weiss \cite{MW03}
but the consistency was  secured only for the  subclass of all stationary and ergodic binary time series with 
almost surely continuous conditional probability $P(X_1=1|\dots,X_{-1},X_0 )$.

However for the class of all stationary and ergodic Markov-chains of some  finite order Problem~1 can be solved.
 Indeed,  
if the time series is a Markov-chain of some finite order, we can estimate the order 
(e.g. as in Csisz\'ar, Shields 
\cite{CsSh00} and Csisz\'ar \cite{Csiszar02}) and  count frequencies of  blocks with length equal to the order.  
Bailey showed that one can't test for being in the class. 

\bigskip
It is conceivable that one can improve the result  of Morvai \cite{Mo00} or Morvai and Weiss \cite{MW03}
so that if the process happens to be Markovian 
then one eventually estimates at all times. Our purpose in this paper is to show that this is not possible. 
This puts some new restrictions on what can be achieved in estimating along stopping times. 

\bigskip
\noindent
\begin{theorem} \label{Thlimits}
 For any strictly increasing sequence of stopping times $\{\lambda_n\}$ 
such that for all stationary and ergodic binary 
Markov-chains with arbitrary finite order,  eventually $\lambda_{n+1}=\lambda_n+1$, 
 and for any sequence of estimators  
$\{h_n(X_{0},\dots,X_{\lambda_n})\}$  
there is a stationary and ergodic  binary time series 
 $\{X_n\}$ with almost surely continuous conditional probability 
$P(X_1=1|\dots,X_{-1},X_0 )$, 
 such that
$$
P\left( \limsup_{n\to\infty} 
| h_n(X_0,\dots,X_{\lambda_n})- P(X_{\lambda_n+1}=1|X_0,\dots,X_{\lambda_n})|>0 \right) 
>0.
$$
\end{theorem}

\bigskip
\noindent
{\bf Remark:} 
\noindent
Bailey \cite{Bailey76} among other things proved that there is no sequence of functions 
$\{e_n(X_0^{n-1})\}$ which for all stationary and ergodic time series, 
if it turns out to be a Markov-chain, would be eventually $1$ and $0$ otherwise. 
(That is, there is no test for the Markov property.)
This result does not imply ours. On the other hand, our result implies Bailey's. 
(Indeed, if there were a test for Markov-chains 
in the above sense, we could apply the estimator in 
Morvai \cite{Mo00} or Morvai and Weiss \cite{MW03} if the time series is not a Markov-chain of some finite order, 
and if the time series is a Markov-chain of some finite order we can estimate the 
order of the Markov chain (e.g. as in Csisz\'ar, Shields 
\cite{CsSh00} or Csisz\'ar \cite{Csiszar02}) and  count frequencies of  blocks with length equal to the order.  

 Bailey \cite{Bailey76} and Ryabko \cite{Ryabko88} proved less than our Theorem~\ref{Thlimits}. 
 They proved the nonexistence of the desired  estimator when the  estimator should work for all 
stationary and ergodic binary 
time series and   when all $\lambda_n=n$, that is, when we always require good prediction.

\section{Proof of Theorem~\ref{Thlimits}}

\bigskip
\noindent
{\sc Proof:}

The proof mainly follows the footsteps of Ryabko \cite{Ryabko88} and Gy\"{o}rfi, Morvai, Yakowitz \cite{GYMY98} 
with alterations where necessary.  For $m\le n$ let $X_m^n=(X_m,\dots,X_n)$.
First  we  define the same  
Markov-chain as in  Ryabko \cite{Ryabko88}  which serves as the technical tool for  construction of our 
counterexample. Let the state space $S$ be the non-negative integers. From state $0$ the process certainly passes to state $1$ and then 
to state $2$, at the following epoch. From each state $s\ge 2$,  the 
Markov chain  passes either to state $0$ or to state 
$s+1$ with equal probabilities $0.5$.  
This construction yields a stationary and ergodic 
Markov chain  $\{M_i\}$ with stationary distribution
$$
P(M=0)=P(M=1)={1\over 4}
$$
and
$$
P(M=i)={1\over  2^{i}} \mbox{\ \  for $i\ge 2$}.
$$

\noindent                          
Let $\psi_k$ denote the first positive time of  occurrence of 
state $2k$ :  
$$ \psi_k=
\min\{i\ge 0: M_i=2k\}.
$$
Note that if $M_0=0$ then $M_i\le 2k$ for $0\le i\le \psi_k$. 
For each $0\le j < \infty$ we will define a binary-valued  Markov-chain  
$\{X^{(j)}_i\}$ with some finite order,  
which we denote as $X^{(j)}_i=f^{(j)}(M_i)$ where $f^{(j)}$ will be a 
$\{0,1\}$ valued function of the state space $S$.  
We will also define a process $\{X_i\}$ which we denote as 
$X_i=f^{(\infty)}(M_i)$ where $f^{(\infty)}$ is also a binary valued function of 
the state space $S$, and the time series $\{X_i\}$  will serve as the stationary 
(non Markov ) unpredictable process. 
For all $0\le j\le \infty$, let $f^{(j)}(0)=0$, $f^{(j)}(1)=0$, and $f^{(j)}(s)=1$ for 
all even states $s$. Note that so far we have only defined $f^{(j)}$ partially. 
We will define the values for the remaining states later on.   
A feature of this definition of $f^{(j)}(\cdot)$ is that whenever  
$X^{(j)}_n=0,X^{(j)}_{n+1}=0,X^{(j)}_{n+2}=1$ we  know 
 that $M_n=0$ and vice versa. 

Now observe that if for a certain  $0\le j\le \infty $, there is an index $K_j$ such that  $f^{(j)}(i)=1$ for all $i\ge K_j$ then the defined process $\{X^{(j)}_n\}$ is a binary Markov-chain with order not greater than $K_j$. 
(Indeed, the probabilities $P(X^{(j)}_n=1|X^{(j)}_0, \dots, X^{(j)}_{n-1})$ are 
determined by the last $K_j$ bits $(X^{(j)}_{n-K_j},\dots,X^{(j)}_{n-1})$. 
To see this consider the following cases. 
\begin{itemize}
\item[I.] If for some $1\leq i\leq K_j-2$ $X^{(j)}_{n-i}=1$ and $X^{(j)}_{n-1-i}=X^{(j)}_{n-2-i}=0$ 
than we can detect that  $M_{n-i}=2$, $M_{n-1-i}=1$ and $M_{n-2-i}=0$  
and the conditional probability does not depend on previous values. 

\item[II.]  If there is no  $1\leq i\leq K_j-2$ such that $X^{(j)}_{n-i}=1$ and 
$X^{(j)}_{n-1-i}=X^{(j)}_{n-2-i}=0$ we have three sub-cases. 
  \begin{itemize} 
     \item[II/1.]  If  $X^{(j)}_{n-1}=1$ then $M_{n-1}\geq K_j$. In this case the 
conditional probability is 
$0.5$.
 
	\item[II/2.] If $X^{(j)}_{n-2}=X^{(j)}_{n-1}=0$ then $M_{n-1}=1$ and the conditional 
	probability is $1$. 

\item[II/3.] If $X^{(j)}_{n-2}=1$ and $X^{(j)}_{n-1}=0$ then $M_{n-1}=0$ 
and so the conditional probability is $0$.)    
	\end{itemize}
\end{itemize}

Now let $f^{(0)}(2k+1)=1$ for all $k\ge 1$ and so the function $f^{(0)}$ is fully defined.  Since $f^{(0)}(i)$ is eventually $1$, the defined process $\{X^{(0)}_i\}$ is a stationary ergodic binary Markov-chain with some finite order.

\noindent 
For function $f^{(j)}$ and index $2k$, if $f^{(j)}(i)$ is defined for all $0\le i\le 2k$, then 
it is easy to see that if  $M_0=0$ (that is, $f^{(j)}(M_0)=0$, $f^{(j)}(M_1)=0$, $f^{(j)}(M_2)=1$ ) 
then $M_i\le 2k$ for $0\le i\le \psi_k$ and 
the mapping 
$$M_0^{\psi_k}\rightarrow(f^{(j)}(M_0),\dots,f^{(j)}(M_{\psi_k}))$$
is invertible. 
If we let $\lambda_n$ operate on process $\{X^{(j)}_i\}$, define 
$$
A_j(k)=\{M_0=0,\psi_{k}=\lambda_n(X^{(j)}_0,X^{(j)}_1,\dots) \ \ \mbox{for some n}\}.
$$
Thus as soon as   $f^{(j)}(i)$ is defined for all $0\le i\le 2k$  
the set $A_j(k)$ is also well defined, it is measurable with respect to $M_0^{\psi_{k}}$ 
and depends on state $2k$ and index $j$ which selects the process $\{X^{(j)}_n\}$ on which 
the stopping times  $\{\lambda_n\}$ operate. 

\noindent
Let $N_{-1}=1$. Notice that $A_0(k)$ is well defined for all $k$. 
Now 
we define $f^{(j)}$ by induction. Assume 
that for $0\le i\le j-1$ 
we have already defined 
a strictly  increasing sequence of integers $N_{i-1}$,  and functions $f^{(i)}$
which are eventually constant. 

\noindent
Now 
we define $f^{(j)}$. Since by assumption $\{X^{(j-1)}_n\}$ 
is a stationary and ergodic binary-valued Markov process with some finite order, 
the estimator is assumed to predict eventually on this process and there is a 
$N_{j-1}> N_{j-2}$ such that 
\[
P(A_{j-1}(N_{j-1}))>1/8.
\] 
Now 
for each $j\le l\le \infty$  
define $f^{(l)}(2m+1)$ for the segment $N_{j-2}\le m< N_{j-1}$ as follows, 
$$
 f^{(l)}(2m+1)= f^{(j-1)}(2m+1). 
$$
Notice that now $A_j(N_{j-1})$ is well defined and coincides with $A_{j-1}(N_{j-1})$. 
We will define $f^{(j)}(2N_{j-1}+1)$ maliciously.   Let
$$
B_j^+= A_j(N_{j-1})\bigcap \{ h_n(f^{(j)}(M_0),\dots,f^{(j)}(M_{\psi_{N_{j-1} }}))\ge {1\over 4} \}
$$
and 
$$
B_j^-= A_j(N_{j-1})\bigcap \{h_n(f^{(j)}(M_0),\dots,f^{(j)}(M_{\psi_{N_{j-1}}}))<  {1\over 4} \}.
$$
 Now notice that the sets $B_j^+$ and $B_j^-$ do not depend on the future values of $f^{(j)}(2r+1)$ for $r\ge N_{j-1}$.  
One of the two sets $B_j^+$, $B_j^-$ has at least probability $1/16$. 
Now we specify $f^{(j)}(2N_{j-1}+1)$. 
Let $f^{(j)}(2N_{j-1}+1)=1$, $I_j=B_j^-$  if $P(B_j^-)\ge P(B_j^+)$ and
let $f^{(j)}(2N_{j-1}+1)=0$, $I_j=B_j^+$ if  $P(B_j^-)< P(B_j^+)$. 
 
\noindent
 Because of the construction
of $\{M_i\}$,  on event $I_j$,
\begin{eqnarray*}
\lefteqn{P(X^{(j)}_{\psi_{N_{j-1}}+1}=1|X^{(j)}_0,\dots, X^{(j)}_{\psi_{N_{j-1}}})}\\
&=&
f^{(j)}(2N_{j-1}+1) P(X^{(j)}_{\psi_{N_{j-1}}+1}=f(2N_{j-1}+1)| X^{(j)}_0, 
\dots,X^{(j)}_{\psi_{N_{j-1}}})\\
&=& f^{(j)}(2N_{j-1}+1) P(M_{\psi_{N_{j-1}+1}}=2N_{j-1}+1| M_0^{\psi_{N_{j-1}}})\\
&=& 0.5 f^{(j)}(2N_{j-1}+1).
\end{eqnarray*}
The difference of the  estimate and the conditional probability is at least ${1\over 4}$ on set $I_j$ and this event occurs  with probability not less  than  $1/16$. 

Now for all $N_{j-1}<m$ define 
$$
f^{(j)}(2m+1)=1.
$$
In this way, $\{X_i^{(j)}\}$ is also a stationary and ergodic binary-valued Markov-chain. 

Now by induction, we defined all the functions $f^{(j)}$ for $0\le j<\infty$. Since $f^{(\infty)}(m)=f^{(j)}(m)=f^{(j-1)}(m)$ for all $0\le m\le 2N_{j-1}$ so we also defined $f^{(\infty)}$. 

Finally by Fatou's Lemma, 
\begin{eqnarray}
\nonumber
\lefteqn{P(\limsup_{n\to\infty} 
\{|h_n(X_0^{\lambda_n})- P(X_{\lambda_n+1}=1|X_0^{\lambda_n})|\ge 1/4\})}\\
\nonumber
&\ge&
P(\limsup_{j\to\infty} I_j)\ge\limsup_{j\to\infty} P(I_j)\ge {1\over 16}.
\end{eqnarray} 
Concerning the conditional probability $P(X_1=1|X_{-\infty}^0)$  observe that as soon as one finds the pattern '001' in 
the sequence $X_{-\infty}^0$ the conditional probability does not depend on previous values. 
The probability of the  occurence of $'001'$ in the past is one since the original Markov 
chain is  ergodic and our process is therefore also ergodic. 
Thus the conditional probabilities are almost surely continuous.   
The proof of Theorem~\ref{Thlimits} is complete.


\end{document}